\documentclass[10pt]{amsart}
\usepackage{enumerate,amsmath,amssymb,latexsym,
amsfonts, amsthm, amscd, mathabx}

\theoremstyle{plain}

\newtheorem{theorem}{Theorem}

\numberwithin{equation}{section}

\addtocounter{section}{-1}

\begin{document}

\title {Closed graphs and open maps}

\date{}

\author [P.L. Robinson]{P.L. Robinson}

\address{Department of Mathematics \\ University of Florida \\ Gainesville FL 32611  USA }

\email[]{paulr@ufl.edu}

 \keywords{}

\begin{abstract}

We offer a new perspective on the closed graph theorem and the open mapping theorem for separated barrelled spaces and fully complete spaces. 

\end{abstract}

\maketitle

\section*{Introduction} 

\medbreak 

The `classical' closed graph theorem asserts that a linear function between two Banach spaces is continuous if its graph is closed; the `classical' open mapping theorem asserts that a surjective continuous linear function between Banach spaces is open. Although these two theorems are very different in form, they are traditionally regarded as being equivalent, since either is an essentially elementary consequence of the other. In [1] these two theorems are generalized to a context in which they become equivalent in a very precise sense, and the two classical theorems become special cases of one underlying (or overarching) theorem. 

\medbreak 

The `classical' theorems themselves have been generalized considerably. Let $X$ be a separated barrelled space and $Y$ a fully complete space: a general version of the closed graph theorem asserts that a linear function from $X$ to $Y$ is continuous if its graph is closed; a general version of the open mapping theorem asserts that a surjective continuous linear function from $Y$ to $X$ is open. For these general versions, see [2] Chapter VI; see also [3] Chapter 12-5 and Chapter 12-4. Our aim in the present paper is to realize these two general versions as special cases of one `master theorem'. 

\medbreak 

As preparation, it is convenient to set out here some notation pertaining to relations. Let $X$ and $Y$ be sets; a relation from $X$ to $Y$ is a subset $\Gamma \subseteq X \times Y$ of their product. When $A \subseteq X$ and $B \subseteq Y$ we write 
$$A \, \Gamma = \{y : (\exists a \in A) \, (a, y) \in \Gamma \} \subseteq Y$$
and 
$$ \Gamma  B= \{x : (\exists b \in B) \, (x, b) \in \Gamma \} \subseteq X.$$
In particular, $X \, \Gamma = {\rm Ran} \, Y$ is the range of $\Gamma$ and $\Gamma Y = {\rm Dom} \, \Gamma$ is the domain of $\Gamma.$ As a special case, $\Gamma$ might be the graph ${\rm gra} \, \phi = \{ (x, \phi (x)) : x \in X \} \subseteq X \times Y$ of a function $\phi : X \to Y$; in this case, $A \, \Gamma$ is the direct image $\phi \, (A)$ and $\Gamma B$ is the inverse image $\phi^{-1} \, (B)$. 
\medbreak 

\section*{Theorem}

\medbreak 

Let $X$ be a separated (or Hausdorff) barrelled space and $Y$ a fully complete (or Ptak) space. 

\medbreak 

\noindent
{\bf Theorem.} 
{\it Let $\Gamma \subseteq X \times Y$ be a closed linear relation with ${\rm Dom}\, \Gamma = X$. If $B \subseteq Y$ is open then $\Gamma \, B \subseteq X$ is open}.  

\medbreak 

This theorem has the closed graph theorem (CGT) and the open mapping theorem (OMT) as immediate corollaries. 

\medbreak 

\noindent
{\bf Closed graph theorem}. {\it If  $\phi : X \to Y$ is a linear function with closed graph, then $\phi$ is continuous.} 

\begin{proof} 
Apply the Theorem when $\Gamma$ is the graph of $\phi$. The condition ${\rm Dom}\, \Gamma = X$ holds by definition of function. Note that here, if $B \subseteq Y$ then $\Gamma \, B = \phi^{-1} \, (B)$. 
\end{proof} 

\medbreak 

\noindent
{\bf Open mapping theorem}. {\it  If  $\psi : Y \to X$ is a surjective linear function with closed graph, then $\psi$ is open.} 

\begin{proof} 
Apply the Theorem when $\Gamma$ is the transpose of the graph of $\psi$. The condition ${\rm Dom}\, \Gamma = X$ holds by definition of surjectivity. Note that here, if $B \subseteq Y$ then $\Gamma \, B = \psi \, (B)$. 
\end{proof} 

\medbreak 

These versions of CGT and OMT appear in Chapter 12 of [3]: CGT as Theorem 12-5-7 and OMT as Theorem 12-4-9. They also appear in Chapter VI of [2]: CGT as Theorem 6 and OMT as Theorem 7; there, Theorem 6 stems from Proposition 10(i) and Theorem 7 stems from Proposition 10(ii); also, Proposition 10(ii) depends on Proposition 10(i). In practical terms, our proof of the main Theorem simplifies the nature and the function of this Proposition 10. The theoretical value of our perspective is perhaps greater than the practical: the viewpoint of the present paper is that both CGT and OMT are simply special cases of one `master theorem'. 

\medbreak 

\section*{Proof} 

\medbreak 

Here we prove the main Theorem, preparing the way with a Lemma and a Proposition.

\medbreak 

We begin quite generally. Let $X$ and $Y$ be vector spaces and let $\Gamma$ be a linear subspace of the product $X \times Y$ such that ${\rm Dom}\,\Gamma = X$. The embedding $i_Y : Y \to X \times Y : y \mapsto (0, y)$ being a linear map, the inverse image $Y_{\Gamma} =  i_Y^{-1} (\Gamma)$ is a linear subspace of $Y$. Let $x \in X$: as $X = {\rm Dom} \, \Gamma$ there exists $y \in Y$ such that $(x, y) \in \Gamma$; if also $y' \in Y$ and $(x, y') \in \Gamma$ then $(0, y' - y) = (x, y') - (x, y) \in \Gamma$ so that $y' - y \in Y_{\Gamma}$. Accordingly, a linear function 
$$\gamma : X \to Y/Y_{\Gamma}$$ 
is well-defined by the rule that if $x \in X$ then 
$$\gamma (x) = y + Y_{\Gamma}$$ 
for any choice of $y \in Y$ such that $(x, y) \in \Gamma$. It is straightforward to check that if $B \subseteq Y$ then 
$$\Gamma \, B = \gamma^{-1} (B + Y_{\Gamma}).$$

\medbreak 

Now, let $X$ and $Y$ be locally convex topological vector spaces. 

\medbreak 

\noindent
{\bf Lemma.} {\it If $\, \Gamma$ is closed then the graph of $\,\gamma$ is closed.} 

\begin{proof} 
Let $(x_0, y_0  +  Y_{\Gamma}) \in X \times Y/\,Y_{\Gamma}$ be a point outside the graph ${\rm gra} \, \gamma$ of $\gamma$; thus, $(x_0, y_0) \in X \times Y$ is a point outside $\Gamma$. Choose a continuous linear functional $\ell \in (X \times Y)'$ such that $\ell |_{\Gamma} = 0$ but $\ell (x_0, y_0) \neq 0$; see [3] Corollary 7-2-12. Let $q : X \times Y \to X \times Y/\, Y_{\Gamma}$ be the (open) quotient map and let $W = (X \times Y) \setminus {\rm Ker} \, \ell \subseteq X \times Y$ be the open set on which $\ell$ is nonzero. The proof ends with the claim that the open neighbourhood $q(W) \subseteq X \times Y/\, Y_{\Gamma}$ of $(x_0, y_0 + Y_{\Gamma})$ is disjoint from ${\rm gra} \, \gamma$. To justify this claim, let $(x, y + Y_{\Gamma}) \in q(W)$: there exists $(x', y') \in W$ such that $(x', y') - (x, y) \in {\rm Ker}\,q = 0 \times Y_{\Gamma} \subseteq \Gamma \subseteq {\rm Ker}\,\ell$ whence $\ell (x, y) = \ell (x', y') \neq 0$ and therefore $(x, y) \notin \Gamma$; it follows that $(x, y + Y_{\Gamma}) \notin {\rm gra} \, \gamma$. 
\end{proof} 

\medbreak 

Next, assume that $X$ and $Y$ are separated locally convex topological vector spaces, and that the linear subspace $\Gamma \subseteq X \times Y$ is closed. Assume further that $Y$ is fully complete. Denote by $\mathcal{N}(X)$ the set of (not necessarily open) neighbourhoods of zero in $X$. 

\medbreak 

\noindent
{\bf Proposition.} {\it  If 
$$B \in \mathcal{N}(Y) \Rightarrow \overline{\Gamma \, B} \in \mathcal{N}(X)$$ 
then 
$$B \subseteq Y \; {\rm open} \; \Rightarrow \; \Gamma \, B \subseteq X \; {\rm open}.$$}

\begin{proof} 
As the subspace $Y_{\Gamma} =  i_Y^{-1} (\Gamma) \subseteq Y$ is closed and $Y$ is fully complete, the quotient $Y/\,Y_{\Gamma}$ is also fully complete: see the Corollary to Proposition 9 in Chapter VI of [2]; see also Theorem 12-4-5 in [3]. Each $\beta \in \mathcal{N}(Y/\,Y_{\Gamma})$ has the form $\beta = B + Y_{\Gamma}$ for some $B \in \mathcal{N}(Y)$ and the identity recorded prior to the Lemma implies that $\overline{\gamma^{-1} (\beta)} = \overline{\Gamma \, B}$; thus $\beta \in \mathcal{N}(Y/\,Y_{\Gamma}) \Rightarrow \overline{\gamma^{-1} (\beta)} \in \mathcal{N}(X)$ and so $\gamma: X \to Y/\,Y_{\Gamma}$ is nearly continuous in the sense of [2]. With the Lemma, all the pieces are now in place for an application of Proposition 10(i) in Chapter VI of [2]: the linear map $\gamma$ is continuous, whence if $B \subseteq Y$ is open then so is $\Gamma \, B = \gamma^{-1} (B + Y_{\Gamma}) \subseteq X$. 
\end{proof} 

\medbreak 

Thus, the two parts of Proposition 10 in Chapter VI of [2] become one when lifted from the context of linear functions to the context of linear relations: they are two special cases of one `master proposition'.  

\medbreak 

Finally, we prove our main Theorem, which we restate for convenience. Recall the context: $X$ is a separated barrelled space and $Y$ is a fully complete space. 

\medbreak 

\noindent 
{\bf Theorem.} 
{\it Let $\Gamma \subseteq X \times Y$ be a closed linear relation with ${\rm Dom}\, \Gamma = X$. If $B \subseteq Y$ is open then $\Gamma \, B \subseteq X$ is open}

\begin{proof} 
Let $B \in \mathcal{N}(Y)$ be absolutely convex, so that $B$ is absolutely convex and absorbing: as $\Gamma \subseteq X \times Y$ is linear and has the whole of $X$ for its domain, it follows that $\Gamma \, B \subseteq X$ is absolutely convex and absorbing, whence its closure $\overline{\Gamma \, B} \subseteq X$ is a barrel; as $X$ is a barrelled space, it follows that $ \overline{\Gamma \, B} \in \mathcal{N}(X)$. As $\mathcal{N}(Y)$ has a base consisting of absolutely convex neighbourhoods of zero, an application of the Proposition concludes the proof of the Theorem. 
\end{proof} 

\medbreak

\bigbreak 

\begin{center} 
{\small R}{\footnotesize EFERENCES}
\end{center} 
\medbreak 

[1] R.S. Monahan and P.L. Robinson, {\it The closed graph theorem is the open mapping theorem}, arXiv 1912.02626 (2019). 
\medbreak 

[2] A.P. Robertson and W.J. Robertson, {\it Topological Vector Spaces}, Cambridge University Press, Second Edition (1973). 

\medbreak 

[3] A. Wilansky, {\it Modern Methods in Topological Vector Spaces}, McGraw-Hill (1978); Dover Publications (2013). 
\medbreak

\end{document}